\documentclass{article}
\usepackage{amsmath, amssymb, graphics}

\newcommand{\mathsym}[1]{{}}
\newcommand{\unicode}{{}}

\begin{document}

\title{Maximum GCD Among Pairs of Random Integers}
\author{R. W. R. Darling $\&$ E. E. Pyle\\
Mathematics Research Group, National Security Agency\\
9800 Savage Road, Fort George G. Meade, Maryland 20755-6515}
\date{\today}
\maketitle

\begin{abstract}\pmb{ ABSTRACT:} Fix \(\alpha >0\), and sample \(N\) integers uniformly at random from \(\left\{1,2,\ldots ,\left\lfloor e^{\alpha
 N}\right\rfloor \right\}\). Given \(\eta >0\), the probability that the maximum of the pairwise GCDs lies between \(N^{2-\eta }\) and \(N^{2+\eta
}\) converges to 1 as \(N\to \infty \). More precise estimates are obtained. This is a Birthday Problem: two of the random integers are likely to
share some prime factor of order \(\left.N^2\right/\log [N]\). The proof generalizes to any arithmetical semigroup where a suitable form of the Prime
Number Theorem is valid.\end{abstract}

\section*{1. Main Result}

Whereas the distribution of the sizes of the prime divisors of a random integer is a well studied subject --- see portions of Billingsley (1999) --- the authors are unaware of any published results on the pairwise Greatest Common Divisors (GCD) among a large collection of random integers.
Theorem 1.1 establishes probabilistic upper and lower bounds for the maximum of these pairwise GCDs.

\subsection*{1.1 Theorem}

\textit{ Suppose} \(\alpha >0\), \textit{ and} \(T_1,\ldots ,T_N\) \textit{ is a random sample, drawn with replacement, from the integers} \(\left\{n\in
\mathbb{N}: n\leq e^{\alpha  N}\right\}\). \textit{ Let} \(\Gamma _{j,k}\) \textit{ denote the Greatest Common Divisor of} \(T_j\)\textit{  and}
\(T_k\). \textit{ For any} \(\eta >0\),

\begin{equation}\lim_{N\to \infty } \mathbb{P}\left[N^{2-\eta }<\underset{1\leq j<k\leq N}{\max }\left\{\Gamma _{j,k}\right\} < N^{2+\eta }\right]=1.\end{equation}
\textit{Indeed there are more precise estimates: for all} \(s\in (0,1)\), \textit{ and} \(b>0\),\textit{  the right side of} (2) \textit{ is finite,
and}

\begin{equation}\mathbb{P}\left[\underset{1\leq j<k\leq N}{\max }\left\{\Gamma _{j,k}\right\}\geq N^{2/s} b^{1/s}\right]\leq  \frac{1}{2 b}\prod
_{p\in \mathcal{P}} \left(1 + \frac{p^s-1}{p^2-p^s}\right),\end{equation}
\textit{where $\mathcal{P}$ denotes the rational primes; while if }\(\Lambda _{j,k}\) \textit{ denotes the largest common prime factor of} \(T_j\)\textit{ and} \(T_k\),\textit{  then for all} { }\(\theta >0\),

\begin{equation}\lim_{N\to \infty } \mathbb{P}\left[\underset{1\leq k<j\leq N}{\max }\left\{\Lambda _{j,k}\right\}<\frac{N^2}{\log \left[N^{\theta
}\right]}\right]\leq e^{-\theta /8 }.\end{equation}

\pmb{ Supplement: }There is an upper bound, similar to (2), for the radical (i.e the largest square-free divisor) \(\text{rad}\left[\Gamma _{j,k}\right]\)
of the GCD:

\begin{equation}\mathbb{P}\left[\underset{1\leq j<k\leq N}{\max }\left\{\text{rad}\left[\Gamma _{j,k}\right]\right\}\geq N^{2/s} b^{1/s}\right]\leq
 \frac{1}{2 b}\prod _{p\in \mathcal{P}} \left(1 -p^{-2}+ p^{s-2}\right).\end{equation}
The proof, which is omitted, uses methods similar to those of Proposition 2.2, based upon a Bernoulli model for occurrence of prime divisors, instead
of a Geometric model for prime divisor multiplicities. For example, when \(s=0.999\), the product on the right side of (4) is approximately 12.44; for the right side of (2), it is approximately 17.64.

\subsection*{1.2 Overview of the Proof of Theorem 1.1}

Let \(Z_i^k\) be a Bernoulli random variable, which takes the value 1 when prime \(p_i\) divides \(T_k\). As a first step towards the proof, imagine
proving a comparable result in the case where \(\left\{Z_i^k,1\leq k\leq N, i\geq 1\right\}\) were independent, and \( \mathbb{P}\left[Z_i^k=1\right]=1\left/p_i\right.\).
The harder parts of the proof arise in dealing with the reality that, for fixed \(k\), \(\left\{Z_i^k, i\geq 1\right\}\) are negatively associated,
and change with \(N\). Convergence of the series

\[\sum _{p\in \mathcal{P}} p^{-2}\log [p]<\infty \]
ensures that the parameter \(\alpha \), which governs the range of integers being sampled, appears neither in (1), (2), nor (3). However the proof
for the lower bound depends crucially on an exponential (in \(N\)) rate of growth in the range, in order to moderate the dependence among \(\left\{Z_i^k,
i\geq 1\right\}\) for fixed \(k\). 

Consider primes as labels on a set of urns; the random variable \(T_j\) contributes a ball to the urn labelled \(p\) if prime \(p\) divides \(T_j\).
The lower bound comes from showing that, with asymptotic probability at least \(1-e^{-\theta /8 }\) , some urn with a label \(p>N^2/\log \left[N^{\theta
}\right]\) contains more than one ball; in that case prime \(p\) is a common divisor of two distinct members of the list \(T_1,\ldots ,T_N\). The
upper bound comes from an exponential moment inequality.

If \(T_1,\ldots ,T_N\) were sampled uniformly without replacement from the integers from 1 to \(N^2\), the lower bound (3) would fail; see the analysis
in Billingsley (1999) of the distribution of the largest prime divisor of a random integer. In the case of sampling from integers from 1 to \(N^r\),
where \(r\geq 3\), the upper bound (2) remains valid, but we do not know whether the lower bound (3) holds or not.

\subsection*{1.3 Generalizations to Arithmetical Semigroups}

Although details will not be given, the techniques used to prove Theorem 1.1 will be valid in the more general context of a commutative semigroup
\(G\) with identity element 1, containing a countably infinite subset \(\mathcal{P}:=\left\{p_1, p_2, \ldots \right\}\) called the \pmb{ primes}
of \(G\), such that every element \(a\neq 1\) of \(G\) has a unique factorization of the form

\[a = \prod _{i\geq 1} p_i^{e_i}, \left(e_1, e_2, \ldots \right) \in  \mathbb{Z}_+^{\infty }\]
where all but finitely many \(\left(e_i\right)\) are zero. Assume in addition that \(G\) is an \pmb{ arithmetical semigroup} in the sense of Knopfmacher
(1990), meaning that there exists a real-valued norm \(|\cdot |\) on \(G\) such that:

$\bullet $ \(|1|=1\), \(\left|p_i\right|>1\) for all \(p_i\in \mathcal{P}\).

$\bullet $ \(|a b| = |a| |b|\) for all \(a,b\in G\).

$\bullet $ The set \(\pi _G[x]:=\left\{i\geq 1: \left|p_i\right|\leq e^x\right\}\) is finite, for each real \(x>0\).

The only analytic condition needed is an abstract form of the Prime Number Theorem (see Knopfmacher (1990), Chapter 6):

\[\lim_{x\to \infty } x e^{-x} \left|\pi _G[x]\right|= 1,\]
used in the proof of Proposition 4.1.This in turn will imply convergence of series such as:

\[\sum _{p\in \mathcal{P}} \log \left[1 + |p|^{s-2}\right], s<1,\]
which appear (in an exponentiated form) in the bound (2). For example, Landau's Prime Ideal Theorem provides such a result in the case where \(G\)
is the set of integral ideals in an algebraic number field, \(\mathcal{P}\) is the set of prime ideals, and \(|a|\) is the norm of \(a\). Knopfmacher
(1990) also studies a more general setting where, for some \(\delta >0\),

\[\lim_{x\to \infty } x e^{-\delta  x} \left|\pi _G[x]\right|= \delta .\]
The authors have not attempted to modify Theorem 1.1 to fit this case.

\section*{2. Pairwise Minima in a Geometric Probability Model}

\subsection*{2.1 Geometric Random Vectors}

Let { }\(\mathcal{P}:=\left\{p_1, p_2, \ldots \right\}\) denote the rational primes \(\{2,3,5,\ldots \}\) in increasing order. Let $\mathcal{I}$
denote the set of non-negative integer vectors \(\left(e_1,e_2,\ldots \right)\) for which \(\sum e_i<\infty \). Let \(X_1,X_2,\ldots \) be (possibly
dependent) positive integer random variables, whose joint law has the property that, for every \(k\in \mathbb{N}\). and every \(\left(e_1,e_2,\ldots
\right)\in \mathcal{I}\) for which \(e_k=0\),

\begin{equation} \mathbb{P}\left[X_k\geq m| \underset{i\neq k}{\cap }\left\{X_i=e_i\right\}\right]\leq \left(\frac{1}{p_k}\right){}^m.\end{equation}

Let \(\zeta \) denote the random vector:

\begin{equation}\zeta :=\left(X_1, X_2, \ldots \right)\in \mathbb{N}^{\mathbb{N}}.\end{equation}

Consider the finite-dimensional projections of \(X_1,X_2,\ldots \) as a general model for prime multiplicities in the prime factorization of a random
integer, without specifying exactly how that integer will be sampled. Let \(\zeta ^{(1)}, \zeta ^{(2)},\ldots ,\zeta ^{(N)}\) be independent random
vectors, all having the same law as \(\zeta \) in (6). Write \(\zeta ^{(k)}\) as\pmb{  }\(\left(X_1^k, X_2^k, \ldots \right)\). Then 

\[L_{j, k}:=\sum _i \min \left\{X_i^k, X_i^j\right\} \log \left[p_i\right]\]
is a model for the log of the GCD of two such random integers. We shall now derive an upper bound for

\[\Delta _N:=\underset{1\leq k<j\leq N}{\max }\left\{L_{j, k}\right\},\]
which models the log maximum of the pairwise GCD among a set of \(N\) $\texttt{"}$large, random$\texttt{"}$ integers.

\subsection*{2.2 Proposition}

\textit{ Assume the joint law of the components of }$\zeta $\textit{  satisfies} (5).

(i) \textit{ For every} \(s\in (0,1)\), \textit{ the following expectation is finite:}

\begin{equation}\mathbb{E}\left[e^{s L_{k,j}}\right] < \prod _i \left(1 + \frac{p_i^s-1}{p_i^2-p_i^s}\right) =:C_s<\infty , s<1.\end{equation}

(ii) \textit{ For any} \(s\in (0,1)\), \textit{ and} \(b>\left.C_s\right/2\), \textit{ for }\textit{ \(C_s\)}\textit{  as in} (7), \textit{ there
is an upper bound:}

\begin{equation}\mathbb{P}\left[\Delta _N\geq \log \left[N^{2/s}\right]+ s^{-1}\log [b]\right]\leq  \frac{C_s}{2 b}<1.\end{equation}

\pmb{ Proof: }Consider first the case where \(X_1,X_2,\ldots \) are independent Geometric random variables, and

\[ \mathbb{P}\left[X_k\geq m\right]=\left(\frac{1}{p_k}\right){}^m, m=1,2,\ldots \]

It is elementary to check that, for \(s\in (0,1)\), and any \(p\in \mathcal{P}\), if \(X^{\prime\prime },X'\) are independent Geometric random variables
with 

\[\mathbb{P}\left[X^{\prime\prime }\geq m\right]=p^{-m}= \mathbb{P}\left[X'\geq m\right], m = 1,2,\ldots ,\]
then their minimum is also a Geometric random variable, which satisfies

\[\mathbb{E}\left[p^{s \min \left\{X^{\prime\prime }, X'\right\}}\right] = 1 + \frac{p^s-1}{p^2-p^s}<1 + p^{s-2}.\]
It follows from the independence assumption that

\[\mathbb{E}\left[e^{s L_{k,j}}\right] = \mathbb{E}\left[\prod _i p_i^{s \min \left\{X_i^k, X_i^j\right\}}\right] = \prod _i \left(1 + \frac{p_i^s-1}{p_i^2-p_i^s}\right)=
C_s.\]
This verifies the assertion (7). Markov's inequality shows that, for any \(s\in (0,1)\)

\[ C_s\geq  e^{s t}\mathbb{P}\left[L_{k,j}\geq t\right].\]
Furthermore

\[\mathbb{P}\left[\underset{1\leq k<j\leq N}{\max }\left\{L_{k,j}\right\}\geq t\right] = \mathbb{P}\left[\underset{1\leq k<j\leq N}{\cup }\left\{L_{k,j}\geq
t\right\}\right]\]

\[\leq  \sum _{1\leq k<j\leq N} \mathbb{P}\left[L_{k,j}\geq t\right] = \frac{N (N-1)}{2}\mathbb{P}\left[L_{k,j}\geq t\right].\]
It follows that, for \(s\in (0,1)\), \(b>0\), and \(t:=s^{-1}\log \left[b N^2\right]\)

\[\mathbb{P}\left[\Delta _N\geq \log \left[N^{2/s}\right]+ s^{-1}\log [b]\right]\leq \frac{N^2}{2}e^{-s t}C_s = \frac{C_s}{2 b}.\]

It remains to consider the case where \(X_1,X_2,\ldots \) satisfies (5), without the independence assumption. Choose a probability space \((\Omega
,\mathcal{F},\mathbb{P})\) on which independent Geometric random variables \(X_1^{\prime },X_2^{\prime },\ldots \) and \(X_1^{\prime\prime },X_2^{\prime\prime
},\ldots \) are defined, such that for all \(i\geq 1\),

\[\mathbb{P}\left[X_i^{\prime\prime }\geq m\right]=p_i^{-m}= \mathbb{P}\left[X_i^{\prime }\geq m\right], m = 1,2,\ldots .\]

We propose to construct \(\zeta ^{(1)}=\left(X_1^1,X_2^1,\ldots \right)\) and \(\zeta ^{(2)}=\left(X_1^2,X_2^2,\ldots \right)\) by induction, on
this probability space \((\Omega ,\mathcal{F},\mathbb{P})\), so that for each \(n\geq 1\), \(\left\{\left(X_i^1,X_i^2\right) 1\leq i\leq n\right\}\)
have the correct joint law, and

\[X_i^1\leq X_i^{\prime }; X_i^2\leq X_i^{\prime\prime },\text{  }i = 1,2,\ldots .\]
Once this is achieved, monotonicity implies

\[\mathbb{E}\left[e^{s L_{1,2}}\right] \leq  \mathbb{E}\left[\prod _i p_i^{s \min \left\{X_i^{\prime }, X_i^{\prime\prime }\right\}}\right],\]
so the desired result will follow from the previous one for independent Geometric random variables.

Since \(\zeta ^{(1)}\) and \(\zeta ^{(2)}\)are independent, it suffices to construct \(\zeta ^{(1)}\) in terms of \(X_1^{\prime },X_2^{\prime },\ldots
\) so that \(X_i^1\leq X_i^{\prime }\) for all \(i\). Let \(\left(U_{i,j}, i\geq 1,j\geq 0\right)\) be independent Uniform\((0,1)\) random variables.
Suppose either \(i=1\), or else some values \(X_1^1=e_1,X_2^1=e_2,\ldots ,X_{i-1}^1=e_{i-1}\) have already been determined. By assumption, there
exists parameters

\[q_{i,k}:= \mathbb{P}\left[X_i\geq k| \underset{j<i}{\cap }\left\{X_j=e_j\right\}\right]\leq \left(\frac{1}{p_i}\right){}^k, k = 1,2,\ldots  .\]

Use these to construct \(X_i^{\prime }\) and \(X_i^1\) as follows:

\[X_i^{\prime }:=\min \left\{k: U_{i,0} U_{i,1} \ldots  U_{i,k}>\left(\frac{1}{p_i}\right){}^k\right\};\]

\[X_i^1:=\min \left\{k: U_{i,0} U_{i,1} \ldots  U_{i,k}>q_{i,k}\right\}\leq X_i^{\prime }.\]

This completes the construction and the proof, giving the result (8). $\square$

\section*{3. Lower Bound for Largest Collision}

\subsection*{3.1 Random Vectors with Independent Components}

Let { }\(\mathcal{P}:=\left\{p_1, p_2, \ldots \right\}\) denote the rational primes \(\{2,3,5,\ldots \}\) in increasing order, and let \(a_j:=\left(\log
\left[p_j\right]\right){}^{1/2}\). Instead of the Geometric model (5), switch to a Bernoulli model in which \(Z_1,Z_2,\ldots \) are independent Bernoulli
random variables, with

\begin{equation}\mathbb{P}\left[Z_j=1\right] := \frac{1}{p_j}.\end{equation}

Let \(\xi \) denote the random vector

\begin{equation}\xi :=\left(a_1Z_1, a_2Z_2, \ldots \right)\in [0,\infty )^{\mathbb{N}}.\end{equation}
under this new assumption, and let \(\xi ^{(1)}, \xi ^{(2)},\ldots ,\xi ^{(N)}\) be independent random vectors, all having the same law as \(\xi
\). Note that \(\xi ^{(1)}\cdot \xi ^{(2)}\) is not a suitable model for the GCD of two random integers, because the independence assumption (9)
is not realistic. However it is a useful context to develop the techniques which will establish the lower bound in Theorem 1.1.

Write\pmb{  }\(\xi ^{(k)}=\left(a_1Z_1^k, a_2Z_2^k, \ldots \right)\).\pmb{  }We seek a lower bound on the log of the largest prime \(p_i\) at which
a $\texttt{"}$collision$\texttt{"}$ occurs, meaning that \(Z_i^j =1=Z_i^k\) for some \(j,k\):

\[\Delta _N^{\prime }:=\underset{1\leq k<j\leq N}{\max }\left\{\underset{i}{\max }\left\{Z_i^j Z_i^k\log \left[p_i\right]\right\} \leq \underset{1\leq
k<j\leq N}{\max }\left\{\xi ^{(k)}\cdot \xi ^{(j)}\right\}.\right.\]

\subsection*{3.2 Proposition}

\textit{ Given} \(\delta \in (0,\infty )\), \textit{ define} \(\varphi _N:=\varphi _N[\delta ]\) \textit{ implicitly by the identity}

\begin{equation}\underset{\varphi _N}{\overset{2 \varphi _N}{\int }}\frac{N^2dx}{2 x^2\log [x]} = \delta .\end{equation}
\textit{Under the assumption of independence of the components of the random vector} (10),

\begin{equation}\lim_{N\to \infty } \mathbb{P}\left[\Delta _N^{\prime } \geq \log \left[\varphi _N[\delta ]\right]\right] \geq 1 - e^{-\delta }.\end{equation}

\pmb{ Remark: }From the integration bounds:

\[\frac{1}{2\varphi _N \log \left[2 \varphi _N\right]}=\frac{1}{\log \left[2 \varphi _N\right]}\underset{\varphi _N}{\overset{2 \varphi _N}{\int
}}\frac{dx}{x^2}<\frac{2 \delta }{N^2}< \frac{1}{\log \left[\varphi _N\right]}\underset{\varphi _N}{\overset{2 \varphi _N}{\int }}\frac{dx}{x^2}=\frac{1}{2\varphi
_N \log \left[\varphi _N\right]}.\]
it follows that \(\varphi _N\), defined in (11), satisfies \(\varphi _N\log \left[\varphi _N\right]/N^2\to 0.25/\delta .\) Hence for all sufficiently
large \(N\), { }\(\varphi _N<\left.N^2\right/2\), and

\begin{equation}\varphi _N>\frac{N^2}{4 \delta  \log \left[2 \varphi _N\right]}>\frac{N^2}{8 \delta  \log [N]}.\end{equation}

The proof uses the following technical Lemma, which the reader may treat as a warm-up exercise for the more difficult Proposition 4.1.

\subsection*{3.3 Lemma}

\textit{ Let} \(\mathcal{P}_N\) \textit{ denote the set of primes }\textit{ \(p\)}\textit{  such that} \(\varphi _N< p\leq 2\varphi _N\). \textit{
Let} \(\left\{Z_p^k, p\in \mathcal{P}_N, 1\leq k\leq N\right\}\) \textit{ be independent Bernoulli random variables, where} \(\mathbb{P}\left[Z_p^k=1\right]
= 1/p\). \textit{ Take} \(D_p:=Z_p^1+\ldots +Z_p^N\). \textit{ Then}

\begin{equation}\lim_{N\to \infty } \mathbb{P}\left[\underset{p\in \mathcal{P}_N}{\cup }\left\{D_p\geq 2\right\}\right] =1 - e^{-\delta  }.\end{equation}

\pmb{ Proof: }Binomial probabilities give:

\[\mathbb{P}\left[D_p\leq 1\right] = \left(1 - \frac{1}{p}\right)^N + \frac{N}{p}\left(1 - \frac{1}{p}\right)^{N-1}=\text{  }\left(1 - \frac{1}{p}\right)^N
\left(1 +\text{  }\frac{N}{p - 1}\right)\]

\[=\text{  }\left(1 - \frac{N}{p} + \frac{N (N-1)}{2 p^2}-\ldots \right) \left(1 +\text{  }\frac{N}{p - 1}\right) = 1 - \frac{N^2}{2 p^2}+ O\left(\frac{N}{\varphi
_N^2}\right) + O\left(\left(\frac{N}{\varphi _N}\right){}^3\right).\]
Independence of \(\left\{Z_p^k, p\in \mathcal{P}_N, 1\leq k\leq N\right\}\) implies independence of \(\left\{D_p, p\in \mathcal{P}_N\right\}\), so

\[\log \left[\mathbb{P}\left[\underset{p\in \mathcal{P}_N}{\cap }\left\{D_p\leq 1\right\}\right]\right] =\text{  }\sum _{p\in \mathcal{P}_N} \log
\left[\mathbb{P}\left[D_p\leq 1\right] \right] \]

\[= \sum _{p\in \mathcal{P}_N} \log \left[1 - \frac{N^2}{2 p^2}\right] + O\left(\frac{N \left|\mathcal{P}_N\right|}{\varphi _N^2}\right) + O\left(\frac{N^3
\left|\mathcal{P}_N\right|}{\varphi _N^3}\right).\]

Using the estimates \(\varphi _N\log \left[\varphi _N\right] = O\left(N^2\right)\), \(\left|\mathcal{P}_N\right|=O\left(\varphi _N/\log \left[\varphi
_N\right]\right)\), and \(p\left/\varphi _N\right.\leq 2\), the last expression becomes

\[= -\sum _{p\in \mathcal{P}_N} \frac{N^2}{2 p^2} +O\left(\frac{N^2}{\varphi _N^2}\right)+ O\left(\frac{N }{\varphi _N\log \left[\varphi _N\right]}\right)
+ O\left(\frac{N^3 }{\varphi _N^2\log \left[\varphi _N\right]}\right).\]

All terms but the first vanish in the limit, while the Prime Number Theorem ensures that

\[\lim_{N\to \infty } \sum _{p\in \mathcal{P}_N} \frac{N^2}{2 p^2}=\delta .\]

Therefore

\[\lim_{N\to \infty }  \mathbb{P}\left[\underset{p\in \mathcal{P}_N}{\cap }\left\{D_p\leq 1\right\}\right] =\text{  }e^{-\delta }.\]

Thus the limit (14) follows. $\square$

\subsubsection*{3.3.1 Proof of Proposition}

According to our model, if \(D_p\geq 2\) for some \(p=p_i\in \mathcal{P}_N\), then there are indices \(1\leq k<j\leq N\) for which \(Z_i^j =1=Z_i^k\).
Since \(\log \left[p_i\right]\geq \log \left[\varphi _N[\delta ]\right]\),

\[\lim_{N\to \infty } \mathbb{P}\left[\Delta _N^{\prime }\geq \log \left[\varphi _N[\delta ]\right]\right]\geq \lim_{N\to \infty } \mathbb{P}\left[\underset{p\in
\mathcal{P}_N}{\cup }\left\{D_p\geq 2\right\}\right] = 1 - e^{-\delta }.\]
This verifies (12). $\square$

\section*{4. Application: Pairwise GCDs of Many Uniform Random Integers}

We shall now prove an analogue of Lemma 3.3 which applies to random integers, dropping the independence { }assumption for the components of the random
vector (10).

\subsection*{4.1 Proposition}

\textit{ Suppose} \(\alpha >0\), \textit{ and} \(T_1,\ldots ,T_N\) \textit{ is a random sample, drawn with replacement, from the integers} \(\left\{n\in
\mathbb{N}: n\leq e^{\alpha  N}\right\}\). \textit{ Given} \(\delta \in (0,\infty )\), \textit{ define} \(\varphi _N:=\varphi _N[\delta ]\) \textit{
implicitly by the identity }(11)\textit{ . Let} \(\mathcal{P}_N\) \textit{ denote the set of primes }\textit{ \(p\)}\textit{  such that} \(\varphi
_N< p\leq 2\varphi _N\); \textit{ for }\textit{ \(p\in \mathcal{P}_N\)}\textit{  let} \(D_p\) \textit{ denote the number of elements of} \(\left\{T_1,\ldots
,T_N\right\}\) \textit{ which are divisible by }\textit{ \(p\)}\textit{ . Then}

\begin{equation}\lim_{N\to \infty } \mathbb{P}\left[\underset{p\in \mathcal{P}_N}{\cup }\left\{D_p\geq 2\right\}\right] =1 - e^{-\delta  }.\end{equation}

\pmb{ Proof: } As noted above, the Prime Number Theorem ensures that

\[\lim_{N\to \infty } \sum _{p\in \mathcal{P}_N} \frac{N^2}{2 p^2}=\delta .\]
More generally, the alternating series for the exponential function ensures that there is an even integer \(d\geq 1\) such that, given \(\epsilon
\in (0,1)\), for all sufficiently large \(N\),

\[1 - e^{-\delta  /(1+\epsilon )}< \sum _{r=1}^d (-1)^{r+1}I_r<1 - e^{-\delta  /(1-\epsilon )}\]
where, for \(\left\{p_1,\ldots ,p_r\right\}\subset \mathcal{P}_N\)

\[I_r\unicode{F39E}:=\sum _{p_1<\ldots  <p_r} \frac{N^{2 r}}{2^r \left(p_1\ldots  p_r\right){}^2}, r = 1, 2, \ldots , d.\]

Because \(\varphi _N/N^2\to 0\), it follows that, for every \(\left\{p_1,\ldots ,p_d\right\}\subset \mathcal{P}_N\),

\[\frac{p_1\ldots  p_d}{e^{\alpha  N}}<\frac{\left(\varphi _N\right){}^d}{e^{\alpha  N}}<e^{2 d \log [N] - \alpha  N}\to 0.\]

Suppose that, for this constant value of \(d\), we fix some \(\left\{p_1,\ldots ,p_d\right\}\subset \mathcal{P}_N\); instead of sampling \(T_1,\ldots
,T_N\) uniformly from integers up to \(e^{\alpha  N}\), sample \(T_1^{\prime },\ldots ,T_N^{\prime }\) uniformly from integers up to

\[p_1\ldots  p_d\left\lfloor e^{\alpha  N}/\left(p_1\ldots  p_d\right)\right\rfloor .\]

>From symmetry considerations, the Bernoulli random variables \(B_1^{\prime },\ldots ,B_d^{\prime }\) are independent, with parameters \(1\left/p_1\right.,\ldots
,1\left/p_d\right.\), respectively where \(B_i^{\prime }\) is the indicator of the event that \(p_i\) divides \(T_1^{\prime }\). By elementary reasoning,

\[\mathbb{P}\left[D_p\geq 2\right] = \frac{N^2}{2 p^2} + O\left(\left(N\left/\varphi _N\right.\right){}^3\right);\]

\[\mathbb{P}\left[D_{p_1}\geq 2, \ldots , D_{p_r}\geq 2\right] = \frac{N^{2 r}}{2^r \left(p_1\ldots  p_r\right){}^2} + O\left(\left(N\left/\varphi
_N\right.\right){}^{2 r +1}\right), r = 1, 2, \ldots , d.\]

If we were to sample \(T_1,\ldots ,T_N\) instead of \(T_1^{\prime },\ldots ,T_N^{\prime }\), the most that such a probability could change is

\[\mathbb{P}\left[\underset{i=1}{\overset{N}{\cup }}\left\{T_i\neq T_i^{\prime }\right\}\right]\leq \frac{N p_1\ldots  p_d}{e^{\alpha  N}}<e^{(2
d+1) \log [N] - \alpha  N}.\]

The same estimate holds for any choice of \(\left\{p_1,\ldots ,p_d\right\}\subset \mathcal{P}_N\). By the inclusion-exclusion formula, taken to the
first \(d\) terms,

\[\mathbb{P}\left[\underset{p\in \mathcal{P}_N}{\cup }\left\{D_p\geq 2\right\}\right] \geq  \sum _{p\in \mathcal{P}_N} \mathbb{P}\left[D_p\geq 2\right]
- \sum _{p_1<p_2} \mathbb{P}\left[D_{p_1}\geq 2,D_{p_2}\geq 2\right] + \ldots  - \sum _{p_1<\ldots  <p_d} \mathbb{P}\left[D_{p_1}\geq 2, \ldots ,
D_{p_d}\geq 2\right]\]

\[= \sum _{r=1}^d (-1)^{r+1}I_r + O\left(\left(N\left/\varphi _N\right.\right){}^3\right) + \left(
\begin{array}{c}
 N \\
 d
\end{array}
\right)e^{(2 d+1) \log [N] - \alpha  N}.\]

So under this simplified model, the reasoning above combines to show that,  for all sufficiently large \(N\),

\[1 - e^{-\delta  /(1+\epsilon )}< \mathbb{P}\left[\underset{p\in \mathcal{P}_N}{\cup }\left\{D_p\geq 2\right\}\right] < 1 - e^{-\delta  /(1-\epsilon
)}.\]

Since $\epsilon $ can be made arbitrarily small, this verifies the result. $\square$

\subsection*{4.2 Proof of Theorem 1.1}

Suppose \(\alpha >0\), and \(T_1,\ldots ,T_N\) is a random sample, drawn with replacement, from the integers \(\left\{n\in \mathbb{N}: n\leq e^{\alpha
 N}\right\}\). Let \(\Lambda _{j,k}\) denote the largest common prime factor of \(T_j\) and \(T_k\). Take

\[\Delta _N^{\prime }:=\underset{1\leq k<j\leq N}{\max }\left\{\log \left[\Lambda _{j,k}\right]\right\}.\]

In the language of Proposition 4.1, if \(D_p\geq 2\) for some \(p\in \mathcal{P}_N\), then there are indices \(1\leq k<j\leq N\) for which \(\Lambda
_{j,k}\text{  }>\varphi _N\). So inequality (13) and Proposition 4.1 imply that, for any \(\theta =8 \delta >0\)

\[\lim_{N\to \infty } \mathbb{P}\left[\Delta _N^{\prime }\geq 2 \log [N] - \log \left[\log \left[N^{\theta }\right]\right]\right]\geq \lim_{N\to
\infty } \mathbb{P}\left[\Delta _N^{\prime }\geq \log \left[\varphi _N[\theta /8]\right]\right]\]

\[\geq \lim_{N\to \infty } \mathbb{P}\left[\underset{p\in \mathcal{P}_N}{\cup }\left\{D_p\geq 2\right\}\right] = 1 - e^{-\theta /8}.\]

This is precisely the lower bound (3). For any \(\eta >0\), the lower bound in (1) follows from:

\[\lim_{N\to \infty } \mathbb{P}\left[\Delta _N^{\prime }>(2-\eta ) \log [N]\right]=1.\]

Let \(\Gamma _{j,k}\geq \Lambda _{j,k}\) denote the Greatest Common Divisor of \(T_j\) and \(T_k\). To obtain the upper bound (2) on \(\Gamma _{j,k}\),
{ }it suffices by Proposition 2.2 to check that condition (5) is valid, when \(X_i\) denotes the multiplicity to which prime \(p_i\) divides \(T_1\).
Take any positive integer \(r\geq 1\), any prime \(p_k\) coprime to \(r\), and any \(m\geq 1\). The conditional probability that \(p_k^m\) divides
\(T_1\), given that \(r\) divides \(T_1\), is

\[\frac{\left\lfloor e^{\alpha  N}/\left(r p_k^m\right)\right\rfloor }{\left\lfloor \left.e^{\alpha  N}\right/r\right\rfloor }\leq  \left(\frac{1}{p_k}\right){}^m.\]

So condition (5) holds. Thus { }(8) holds, which is equivalent to (2). 

Finally we derive the upper bound in (1), for an arbitrary \(\eta >0\). Fix \(\epsilon \in (0,1)\) and \(\eta >0\). Select \(s\in (0,1)\) to satisfy
\(2/s=2 + \eta /2\). Then choose \(b=\left.C_s^{\prime }\right/\epsilon \). According to (8),

\[\mathbb{P}\left[\Delta _N\geq (2 + \eta /2)\log [N]+ s^{-1}\log [b]\right]\leq \epsilon /2.\]

For any \(N\) sufficiently large so that \((\eta /2)\log [N]>s^{-1}\log [b]\), 

\[\mathbb{P}\left[\Delta _N\geq (2 + \eta )\log [N]\right]\leq \epsilon /2.\]

This yields the desired bound (1). $\square$

\end{document}